# Computational procedures for weighted projective spaces

Michele Rossi
Università degli Studi di Torino
Italy
michele.rossi@unito.it

LeaTerracini
Università degli Studi di Torino
Italy
lea.terracini@unito.it

## ▼ Introduction

Let us introduce the main geometric objects the following procedures are dealing with. For any further detail the interested reader is referred to the extensive mathematical treatment of the subject [RT] and to references therein.

A weighted projective space (wps) is a suitable finite quotient $\boldsymbol{P}(Q) := \boldsymbol{P}^n/W_Q$ of the $n$-dimensional projective space, where $Q = (q_0, \ldots, q_n)$ is the *weights vector,* i.e. a $(n+1)$-vector with coprime positive integer entries, characterizing the action defining the quotient [RT, 1.4]. It is probably the better known example of a *singular, complete, normal toric variety*.

Since normal toric varieties can be defined starting from convex geometric objects like cones, fans of cones and, for projective ones, polytopes [RT, 1.1], the approach to wps's via toric geometry naturally provides an algorithmic treatment of the subject. Moreover polytopes associated with polarized wps's turns out to be simplices, meaning that their treatment can be further reduced to consider suitable associated matrices.

Here we will provide procedures which are able to produce the toric data associated with a (polarized) wps i.e. fans , polytopes and their equivalences. But more originally we will provide procedures which are able to detect a weights vector $Q$ starting from either a fan or a polytope: we will call this process *the recognition of a (polarized) wps*. Moreover we will give procedures connecting polytopes of a polarized wps with an associated fan and viceversa.

Although there already exist many mathematical packages which are able to compute the toric data of a given toric variety, and in particular of a given wps (the interested reader is referred to the well updated D. Cox web-page [Cox] and references therein), the recognition of toric data, the connection between fans and polytopes and the detection of equivalence between toric data with the computation of the linking transformations are all original procedures which do not find comparisons elsewhere.

> *restart* :
> *with*(*LinearAlgebra*) :

## ▼ 0.a Preliminaries

# 1. Reduction of weights vector

Given a weights vector Q, computes its reduction (see [RT, 1.6]):

> $ReducedWeights := \mathbf{proc}(Q)$ **local** $nn, dd, i, aa, L : nn := nops(Q) - 1 : dd := [\,] : $ **for** $i$ **from** $1$ **to** $nn + 1$ **do** $dd := [op(dd), igcd(seq(Q[j], j = 1 .. i-1), seq(Q[j], j = i + 1 .. nn + 1))] $ **end do**; $aa := [\,] : $ **for** $i$ **from** $1$ **to** $nn + 1$ **do** $aa := [op(aa), ilcm(seq(dd[j], j = 1 .. i - 1), seq(dd[j], j = i + 1 .. nn + 1))] $ **end do**: $L := [\,] :$ **for** $i$ **from** $1$ **to** $nn + 1$ **do**
> $L := \left[op(L), \dfrac{Q[i]}{aa[i]}\right]$ **end do**: $L$ **end proc**:

**Example:**

> $Q := [3, 2, 4] :$
> $ReducedWeights(Q)$

$$[3, 1, 2] \qquad (2.1.1)$$

# 2. From fan and polytopes to matrices and conversely

The following procedures associate a matrix to the (ordered) 1-skeleton of a fan and viceversa. The same is done for vertices of a polytope.

> $FanMatrix := \mathbf{proc}(fan)\ Transpose(Matrix(fan))\ \mathbf{end\ proc}:$
> $MatrixFan := \mathbf{proc}(V)\ convert(Transpose(V), listlist)\ \mathbf{end\ proc}:$
> $PolytopeMatrix := \mathbf{proc}(P)\ Transpose(Matrix([seq(P[i] - P[1], i = 2 .. nops(P[1]) + 1)]))\ \mathbf{end\ proc}:$
> $MatrixPolytope := \mathbf{proc}(W)\ [[seq(0, i = 1 .. op(W)[1])], op(convert(Transpose(W), listlist))]\ \mathbf{end\ proc}:$

**Example:**

> $fan := [[3, 2, 1], [7, 3, 1], [0, 0, 2], [3, 1, 3]]$

$$fan := [[3, 2, 1], [7, 3, 1], [0, 0, 2], [3, 1, 3]] \qquad (2.2.1)$$

> $FM := FanMatrix(fan)$

$$FM := \begin{bmatrix} 3 & 7 & 0 & 3 \\ 2 & 3 & 0 & 1 \\ 1 & 1 & 2 & 3 \end{bmatrix} \qquad (2.2.2)$$

> $MatrixFan(FM)$

$$[[3, 2, 1], [7, 3, 1], [0, 0, 2], [3, 1, 3]] \qquad (2.2.3)$$

> $P := [[1, 2, 3, 4], [3, 2, 6, 1], [2, 1, 3, 5], [2, 4, 1, 1], [7, 8, 8, 8]]$

$$P := [[1, 2, 3, 4], [3, 2, 6, 1], [2, 1, 3, 5], [2, 4, 1, 1], [7, 8, 8, 8]] \qquad (2.2.4)$$

> $PM := PolytopeMatrix(P)$

$$PM := \begin{bmatrix} 2 & 1 & 1 & 6 \\ 0 & -1 & 2 & 6 \\ 3 & 0 & -2 & 5 \\ -3 & 1 & -3 & 4 \end{bmatrix} \qquad (2.2.5)$$

> $MatrixPolytope(PM)$

$$[[0, 0, 0, 0], [2, 0, 3, -3], [1, -1, 0, 1], [1, 2, -2, -3], [6, 6, 5, 4]] \qquad (2.2.6)$$

## 3. Is Gorenstein i.e. Fano

A wps is Gorenstein if and only if it is Fano. Therefore the following procedure are essentially the same. *IsFano* gives also the index of $P(Q)$ when it is Fano. First of all let us compute *the sum of weights* in $Q$

> $Qsum := \mathbf{proc}(Q)$ **local** $QS, i : QS := 0 :$ **for** $i$ **to** $nops(Q)$ **do** $QS := QS + Q[i]$ **end do**:
$QS$ **end proc**:

> $Q$

$$[3, 2, 4] \quad (2.3.1)$$

> $Qsum(Q)$

$$9 \quad (2.3.2)$$

Then let us compute *the least common multiple of weights* in $Q$

> $Qdelta := \mathbf{proc}(Q) \ lcm(seq(Q[i], i = 1..nops(Q))) \ \mathbf{end \ proc}$

$$Qdelta := \mathbf{proc}(Q) \ lcm(seq(Q[i], i = 1..nops(Q))) \ \mathbf{end \ proc} \quad (2.3.3)$$

> $Qdelta(Q)$

$$12 \quad (2.3.4)$$

Finally the procedures *IsGorenstein* and *IsFano*

> $IsGorenstein := \mathbf{proc}(Q) \ \mathbf{if} \ type\left(\dfrac{Qsum(Q)}{Qdelta(Q)}, integer\right) = true \ \mathbf{then} \ true \ \mathbf{else} \ false \ \mathbf{end \ if}$
**end proc**

$IsGorenstein := \mathbf{proc}(Q)$
  **if** $type(Qsum(Q)/Qdelta(Q), integer) = true$ **then** $true$ **else** $false$ **end if**
**end proc** $\quad (2.3.5)$

> $IsGorenstein(Q)$

$$false \quad (2.3.6)$$

> $IsGorenstein([1, 1, 1, 1])$

$$true \quad (2.3.7)$$

> $IsFano := \mathbf{proc}(Q) \ print\left(IsGorenstein(Q), index = \dfrac{Qsum(Q)}{Qdelta(Q)}\right) \ \mathbf{end \ proc}$

$IsFano := \mathbf{proc}(Q) \ print(IsGorenstein(Q), index = Qsum(Q)/Qdelta(Q)) \ \mathbf{end \ proc} \quad (2.3.8)$

> $IsFano(Q)$

$$false, index = \dfrac{3}{4} \quad (2.3.9)$$

> $IsFano([1, 1, 1, 1])$

$$true, index = 4 \quad (2.3.10)$$

# 0.b Given a weights vector Q, determine a fan associated to Q

Given a weights vector Q, the following procedure computes a fan associated to Q:

> $Fan := \mathbf{proc}(Q)$
**local** $n, U1, U2, SS, TT, VV :$
$n := nops(Q) - 1 :$
$U1 := LinearAlgebra[HermiteForm](Matrix([seq([Q[i]], i = 1..n+1)]), output = ['U'])$ :
$U2 := Matrix([seq(convert(Row(U1, i), list), i = 2..n+1)]) :$
$VV := convert(Transpose(U2), listlist) :$
$VV$

**end proc**:

**Examples:**

> $Q := [3, 2, 4]$

$$Q := [3, 2, 4] \tag{3.1}$$

> $Fan(Q)$

$$[[-2, -2], [3, 1], [0, 1]] \tag{3.2}$$

> $RQ := ReducedWeights(Q)$

$$RQ := [3, 1, 2] \tag{3.3}$$

> $Fan(RQ)$

$$[[-1, -1], [3, 1], [0, 1]] \tag{3.4}$$

> $Q := [1, 1, 2, 4]$

$$Q := [1, 1, 2, 4] \tag{3.5}$$

> $Fan(Q)$

$$[[-1, -2, -4], [1, 0, 0], [0, 1, 0], [0, 0, 1]] \tag{3.6}$$

> $Q := [16, 24, 54, 27, 63, 32, 45, 56, 34, 84, 178, 236, 142, 266, 988, 1016]$

$$Q := [16, 24, 54, 27, 63, 32, 45, 56, 34, 84, 178, 236, 142, 266, 988, 1016] \tag{3.7}$$

> $RedQ := ReducedWeights(Q)$

$$RedQ := [16, 24, 54, 27, 63, 32, 45, 56, 34, 84, 178, 236, 142, 266, 988, 1016] \tag{3.8}$$

> $Fan(Q)$

$$[[-51, -48, -24, -39, -2, -6, -29, -34, -18, -43, -53, -28, -23, -100, -89], [16, 14, 7, 11, 0, 1, 8, 10, 4, 10, 12, 6, 2, 12, 8], [0, 1, 0, 0, 0, 0, 0, 0, 0, 0, 0, 0, 0, 0, 0], [16, 14, 8, 11, 0, 1, 8, 10, 4, 10, 12, 6, 2, 12, 8], [0, 0, 0, 1, 0, 0, 0, 0, 0, 0, 0, 0, 0, 0, 0], [0, 0, 0, 0, 1, 0, 0, 0, 0, 0, 0, 0, 0, 0, 0], [0, 0, 0, 0, 0, 1, 0, 0, 0, 0, 0, 0, 0, 0, 0], [0, 0, 0, 0, 0, 0, 1, 0, 0, 0, 0, 0, 0, 0, 0], [0, 0, 0, 0, 0, 0, 0, 1, 0, 0, 0, 0, 0, 0, 0], [0, 0, 0, 0, 0, 0, 0, 0, 1, 0, 0, 0, 0, 0, 0], [0, 0, 0, 0, 0, 0, 0, 0, 0, 1, 0, 0, 0, 0, 0], [0, 0, 0, 0, 0, 0, 0, 0, 0, 0, 1, 0, 0, 0, 0], [0, 0, 0, 0, 0, 0, 0, 0, 0, 0, 0, 1, 0, 0, 0], [0, 0, 0, 0, 0, 0, 0, 0, 0, 0, 0, 0, 1, 0, 0], [0, 0, 0, 0, 0, 0, 0, 0, 0, 0, 0, 0, 0, 1, 0], [0, 0, 0, 0, 0, 0, 0, 0, 0, 0, 0, 0, 0, 0, 1]] \tag{3.9}$$

## ▼ Given a weights vector Q, determine the Q-*canonical* fan of P(Q)

Given a weights vector Q, *QcanonicalFan* returns the *Q*-canonical fan of $P(Q)$, as defined in [RT, 2.4].

> $QcanonicalFan := \mathbf{proc}(Q)$
> **local** $n, U1, U2, U3, SS, TT, VV, VV1, Q1, UU4, U4, VVV$:
> $n := nops(Q) - 1$:
> $Q1 := [seq(Q[i], i = 2..n + 1), Q[1]]$:
> $U1 := LinearAlgebra[HermiteForm](Matrix([seq([Q1[i]], i = 1..n + 1)]), output = ['U'])$:
> $U2 := Matrix([seq(convert(Row(U1, i), list), i = 2..n + 1)])$:
> $U3 := LinearAlgebra[HermiteForm](U2, output = ['H'])$:
> $VV := convert(Transpose(U3), listlist)$:
> $VV1 := [VV[n + 1], seq(VV[i], i = 1...n)]$:
> $UU4 := Matrix[column](VV1)$:
> $U4 := Transpose(UU4)$:
> $VVV := convert(Transpose(U4), listlist)$:
> $VVV$
> **end proc**:

**Example:**

> `Fan(Q)`

$[[-51, -48, -24, -39, -2, -6, -29, -34, -18, -43, -53, -28, -23, -100, -89], [16,$     **(3.1.1)**
$14, 7, 11, 0, 1, 8, 10, 4, 10, 12, 6, 2, 12, 8], [0, 1, 0, 0, 0, 0, 0, 0, 0, 0, 0, 0, 0, 0, 0],$
$[16, 14, 8, 11, 0, 1, 8, 10, 4, 10, 12, 6, 2, 12, 8], [0, 0, 0, 1, 0, 0, 0, 0, 0, 0, 0, 0, 0, 0,$
$0], [0, 0, 0, 0, 1, 0, 0, 0, 0, 0, 0, 0, 0, 0, 0], [0, 0, 0, 0, 0, 1, 0, 0, 0, 0, 0, 0, 0, 0, 0],$
$[0, 0, 0, 0, 0, 0, 1, 0, 0, 0, 0, 0, 0, 0, 0], [0, 0, 0, 0, 0, 0, 0, 1, 0, 0, 0, 0, 0, 0, 0], [0, 0,$
$0, 0, 0, 0, 0, 0, 1, 0, 0, 0, 0, 0, 0], [0, 0, 0, 0, 0, 0, 0, 0, 0, 1, 0, 0, 0, 0, 0], [0, 0, 0, 0,$
$0, 0, 0, 0, 0, 0, 1, 0, 0, 0, 0], [0, 0, 0, 0, 0, 0, 0, 0, 0, 0, 0, 1, 0, 0, 0], [0, 0, 0, 0, 0, 0,$
$0, 0, 0, 0, 0, 0, 1, 0, 0], [0, 0, 0, 0, 0, 0, 0, 0, 0, 0, 0, 0, 0, 1, 0], [0, 0, 0, 0, 0, 0, 0, 0,$
$0, 0, 0, 0, 0, 0, 1]]$

> $MF := QcanonicalFan(Q)$

$MF := [[-65, -20, -68, -132, -2, -84, -67, -144, -67, -153, -140, -89, -95, -187,$     **(3.1.2)**
$-127], [1, 0, 0, 0, 0, 0, 0, 0, 0, 0, 0, 0, 0, 0, 0], [0, 1, 0, 0, 0, 0, 0, 0, 0, 0, 0, 0, 0, 0,$
$0], [0, 0, 1, 0, 0, 0, 0, 0, 0, 0, 0, 0, 0, 0, 0], [0, 0, 0, 1, 0, 0, 0, 0, 0, 0, 0, 0, 0, 0, 0],$
$[0, 0, 0, 0, 1, 0, 0, 0, 0, 0, 0, 0, 0, 0, 0], [0, 0, 1, 1, 0, 2, 0, 0, 0, 0, 0, 0, 0, 0, 0], [0, 0,$
$0, 0, 0, 0, 1, 0, 0, 0, 0, 0, 0, 0, 0], [0, 0, 0, 0, 0, 0, 0, 1, 0, 0, 0, 0, 0, 0, 0], [0, 0, 0, 0,$
$0, 0, 0, 0, 1, 0, 0, 0, 0, 0, 0], [0, 0, 0, 0, 0, 0, 0, 0, 0, 1, 0, 0, 0, 0, 0], [0, 0, 0, 0, 0, 0,$
$0, 0, 0, 0, 1, 0, 0, 0, 0], [0, 0, 0, 0, 0, 0, 0, 0, 0, 0, 0, 1, 0, 0, 0], [0, 1, 0, 0, 0, 1, 0, 1,$
$0, 1, 0, 1, 2, 0, 0], [0, 0, 0, 1, 0, 1, 0, 1, 1, 1, 1, 0, 1, 2, 0], [1, 0, 1, 1, 0, 0, 1, 1, 0, 1,$
$1, 1, 0, 1, 2]]$

# 0.c Given a weights vector Q and a polarization O(m), determine a polytope associated to (Q,m)

## 1. Transverse and weighted trasverse of a matrix

The *weighted transversion* is a process giving a matrix associated with the polytope of a minimally polarized wps when starting from a given fan of it [RT, 3.1]. Here this is realized by the procedure *WeightedTransverse*.

Let us first introduce some technical procedures. In the procedure *Matri*, $M$ is a $n \times (n+1)$ matrix, $j$ an index with $1 \leq j \leq n+1$ and *Matri(M,j)* is the matrix obtained by eliminating from $M$ the $j$-th column.

> $Matri := \mathbf{proc}(M, j) \ \mathbf{local} \ n, Mi :$
>    $n := op(M)[1] : \ Mi := Minor(Matrix([[seq(1, k = 1..n + 1)], seq(convert(Row(M,$
>      $i), list), i = 1..n)]), 1, j, output = ['matrix'])$
>    **end proc**:

**Example:**

> $M := Matrix([[1, 2, 5, 4], [3, 4, 2, 1], [5, 6, 7, 2]])$

$$M := \begin{bmatrix} 1 & 2 & 5 & 4 \\ 3 & 4 & 2 & 1 \\ 5 & 6 & 7 & 2 \end{bmatrix} \quad \textbf{(4.1.1)}$$

> $Matri(M, 2);$

$$\begin{bmatrix} 1 & 5 & 4 \\ 3 & 2 & 1 \\ 5 & 7 & 2 \end{bmatrix} \quad (4.1.2)$$

The procedure *Vzero* eliminates the first column from a $n \times (n+1)$ matrix by using *Matri*

> $Vzero := \mathbf{proc}(M) \; Matri(M, 1) \; \mathbf{end \; proc}$:

**Example:**

> $Vzero(M)$

$$\begin{bmatrix} 2 & 5 & 4 \\ 4 & 2 & 1 \\ 6 & 7 & 2 \end{bmatrix} \quad (4.1.3)$$

The procedure Transverse computes the *transverse* of a square matrix M, as defined in [RT,1.3]

> $Transverse := \mathbf{proc}(M) \; (LinearAlgebra[Transpose](M))^{-1} \; \mathbf{end \; proc}$:

Finally the procedure *WeightedTransverse* computing the weighted transverse matrix strating from a fan V.

> $WeightedTransverse := \mathbf{proc}(V) \; \mathbf{local} \; n, M, SS, QQ, de : n := nops(V) - 1 : M$
> $\quad := Transpose(Matrix(V)) :$
> $SS := [seq(Determinant(Matri(M, i)), i = 1 .. n + 1)] :$
> $QQ := [seq(\mathrm{abs}(SS[i]), i = 1 .. n + 1)] : de := ilcm(op(QQ)) :$
> $Transpose(MatrixInverse(Vzero(M))) . Matrix\left(1 .. n, 1 .. n, de \cdot Vector\left(\left[seq\left(\frac{1}{QQ[j]}, j\right.\right.\right.\right.$
> $\left.\left.\left.\left. = 2 .. n + 1\right)\right]\right), shape = diagonal\right) \; \mathbf{end \; proc}$:

**Example:**

> $WeightedTransverse([[-2, -2], [3, 1], [0, 1]])$

$$\begin{bmatrix} 2 & -1 \\ 0 & 3 \end{bmatrix} \quad (4.1.4)$$

## 2. Computing the polytope matrix

*QpolMat(Q,m)* computes a polytope matrix associated with $(P(Q), \mathcal{O}(m))$, as defined in [RT, Thm.3.3], starting from a weights vector *Q* and a polarization *m*: if not expressed as an input, the value *m* = 1 is taken by default. *FpolMat* computes the same polytope matrix when starting from a fan.

> $QpolMat := \mathbf{proc}(Q, m :: integer := 1) \; \mathbf{local} \; VQ : VQ := Fan(Q) : m$
> $\quad \cdot WeightedTransverse(VQ) \; \mathbf{end \; proc}$:
> $FpolMat := \mathbf{proc}(V, m :: integer := 1) \; m \cdot WeightedTransverse(V) \; \mathbf{end \; proc}$:

**Examples:**

> $Q$

$$[16, 24, 54, 27, 63, 32, 45, 56, 34, 84, 178, 236, 142, 266, 988, 1016] \quad (4.2.1)$$

> $F := Fan(Q)$

$F := [[-51, -48, -24, -39, -2, -6, -29, -34, -18, -43, -53, -28, -23, -100, -89],$ (4.2.2)
$\quad [16, 14, 7, 11, 0, 1, 8, 10, 4, 10, 12, 6, 2, 12, 8], [0, 1, 0, 0, 0, 0, 0, 0, 0, 0, 0, 0, 0,$

0], [16, 14, 8, 11, 0, 1, 8, 10, 4, 10, 12, 6, 2, 12, 8], [0, 0, 0, 1, 0, 0, 0, 0, 0, 0, 0, 0, 0, 0, 0], [0, 0, 0, 0, 1, 0, 0, 0, 0, 0, 0, 0, 0, 0, 0], [0, 0, 0, 0, 0, 1, 0, 0, 0, 0, 0, 0, 0, 0, 0], [0, 0, 0, 0, 0, 0, 1, 0, 0, 0, 0, 0, 0, 0, 0], [0, 0, 0, 0, 0, 0, 0, 1, 0, 0, 0, 0, 0, 0, 0], [0, 0, 0, 0, 0, 0, 0, 0, 1, 0, 0, 0, 0, 0, 0], [0, 0, 0, 0, 0, 0, 0, 0, 0, 1, 0, 0, 0, 0, 0], [0, 0, 0, 0, 0, 0, 0, 0, 0, 0, 1, 0, 0, 0, 0], [0, 0, 0, 0, 0, 0, 0, 0, 0, 0, 0, 1, 0, 0, 0], [0, 0, 0, 0, 0, 0, 0, 0, 0, 0, 0, 0, 1, 0, 0], [0, 0, 0, 0, 0, 0, 0, 0, 0, 0, 0, 0, 0, 1, 0], [0, 0, 0, 0, 0, 0, 0, 0, 0, 0, 0, 0, 0, 0, 1]]

> $QpolMat(Q)$

$$\begin{bmatrix} \text{15 x 15 Matrix} \\ \text{Data Type: anything} \\ \text{Storage: rectangular} \\ \text{Order: Fortran\_order} \end{bmatrix}$$ (4.2.3)

> $FpolMat(F, 1)$

$$\begin{bmatrix} \text{15 x 15 Matrix} \\ \text{Data Type: anything} \\ \text{Storage: rectangular} \\ \text{Order: Fortran\_order} \end{bmatrix}$$ (4.2.4)

> $QpolMat(ReducedWeights(Q), 3)$

$$\begin{bmatrix} \text{15 x 15 Matrix} \\ \text{Data Type: anything} \\ \text{Storage: rectangular} \\ \text{Order: Fortran\_order} \end{bmatrix}$$ (4.2.5)

## 3. Computing the polytope

*Qpolytope(Q,m)* computes a polytope defining $(P(Q), \mathcal{O}(m))$; *Fpolytope* computes the same polytope starting from a fan of $P(Q)$ and a polarization $m$. Again the polarization is considered equal to 1 when not expressed as an input.

> $Qpolytope := \textbf{proc}(Q, m :: integer := 1) \ [[seq(0, i = 1..nops(Q) - 1)], \\ op(convert(Transpose(QpolMat(Q, m)), listlist))\ ]$
>  $\textbf{end proc}:$

> $Fpolytope := \textbf{proc}(fan, m :: integer := 1) \ [[seq(0, i = 1.. nops(fan[1]))], \\ op(convert(Transpose(FpolMat(fan, m)), listlist))\ ] \ \textbf{end proc}:$

**Examples:**

> $Qpolytope([3, 2, 4])$

$$[[0, 0], [2, 0], [-1, 3]]$$ (4.3.1)

> $Q$

$$[16, 24, 54, 27, 63, 32, 45, 56, 34, 84, 178, 236, 142, 266, 988, 1016]$$ (4.3.2)

> $F := Fan(Q)$

$F := [[-51, -48, -24, -39, -2, -6, -29, -34, -18, -43, -53, -28, -23, -100, -89],$ (4.3.3)
  $[16, 14, 7, 11, 0, 1, 8, 10, 4, 10, 12, 6, 2, 12, 8], [0, 1, 0, 0, 0, 0, 0, 0, 0, 0, 0, 0, 0, 0,$

0], [16, 14, 8, 11, 0, 1, 8, 10, 4, 10, 12, 6, 2, 12, 8], [0, 0, 0, 1, 0, 0, 0, 0, 0, 0, 0, 0, 0, 0, 0], [0, 0, 0, 0, 1, 0, 0, 0, 0, 0, 0, 0, 0, 0, 0], [0, 0, 0, 0, 0, 1, 0, 0, 0, 0, 0, 0, 0, 0, 0], [0, 0, 0, 0, 0, 0, 1, 0, 0, 0, 0, 0, 0, 0, 0], [0, 0, 0, 0, 0, 0, 0, 1, 0, 0, 0, 0, 0, 0, 0], [0, 0, 0, 0, 0, 0, 0, 0, 1, 0, 0, 0, 0, 0, 0], [0, 0, 0, 0, 0, 0, 0, 0, 0, 1, 0, 0, 0, 0, 0], [0, 0, 0, 0, 0, 0, 0, 0, 0, 0, 1, 0, 0, 0, 0], [0, 0, 0, 0, 0, 0, 0, 0, 0, 0, 0, 1, 0, 0, 0], [0, 0, 0, 0, 0, 0, 0, 0, 0, 0, 0, 0, 1, 0, 0], [0, 0, 0, 0, 0, 0, 0, 0, 0, 0, 0, 0, 0, 1, 0], [0, 0, 0, 0, 0, 0, 0, 0, 0, 0, 0, 0, 0, 0, 1]]

> $Qpolytope(Q)$

[[0, 0, 0, 0, 0, 0, 0, 0, 0, 0, 0, 0, 0, 0, 0], [125253685073790, 0, −250507370147580, 0, 0, 0, 0, 0, 0, 0, 0, 0, 0, 0, 0], [−97419532835170, 111336608954480, 0, 0, 0, 0, 0, 0, 0, 0, 0, 0, 0, 0, 0], [−97419532835170, 0, 222673217908960, 0, 0, 0, 0, 0, 0, 0, 0, 0, 0, 0, 0], [−65609073133890, 0, 0, 95431379103840, 0, 0, 0, 0, 0, 0, 0, 0, 0, 0, 0], [0, 0, 0, 0, 187880527610685, 0, 0, 0, 0, 0, 0, 0, 0, 0, 0], [−8350245671586, 0, 0, 0, 0, 133603930745376, 0, 0, 0, 0, 0, 0, 0, 0, 0], [−53680150745910, 0, 0, 0, 0, 0, 107360301491820, 0, 0, 0, 0, 0, 0, 0, 0], [−110517957418050, 0, 0, 0, 0, 0, 0, 176828731868880, 0, 0, 0, 0, 0, 0, 0], [−17893383581970, 0, 0, 0, 0, 0, 0, 0, 71573534327880, 0, 0, 0, 0, 0, 0], [−21110171641650, 0, 0, 0, 0, 0, 0, 0, 0, 33776274626640, 0, 0, 0, 0, 0], [−19106494333290, 0, 0, 0, 0, 0, 0, 0, 0, 0, 25475325777720, 0, 0, 0, 0], [−15877227685410, 0, 0, 0, 0, 0, 0, 0, 0, 0, 0, 42339273827760, 0, 0, 0], [−2825271091890, 0, 0, 0, 0, 0, 0, 0, 0, 0, 0, 0, 22602168735120, 0, 0], [−4563899456130, 0, 0, 0, 0, 0, 0, 0, 0, 0, 0, 0, 0, 6085199274840, 0], [−2958748466310, 0, 0, 0, 0, 0, 0, 0, 0, 0, 0, 0, 0, 0, 5917496932620]]  **(4.3.4)**

> $Fpolytope(F, 4)$

[[0, 0, 0, 0, 0, 0, 0, 0, 0, 0, 0, 0, 0, 0, 0], [501014740295160, 0, −1002029480590320, 0, 0, 0, 0, 0, 0, 0, 0, 0, 0, 0, 0], [−389678131340680, 445346435817920, 0, 0, 0, 0, 0, 0, 0, 0, 0, 0, 0, 0, 0], [−389678131340680, 0, 890692871635840, 0, 0, 0, 0, 0, 0, 0, 0, 0, 0, 0, 0], [−262436292535560, 0, 0, 381725516415360, 0, 0, 0, 0, 0, 0, 0, 0, 0, 0, 0], [0, 0, 0, 0, 751522110442740, 0, 0, 0, 0, 0, 0, 0, 0, 0, 0], [−33400982686344, 0, 0, 0, 0, 534415722981504, 0, 0, 0, 0, 0, 0, 0, 0, 0], [−214720602983640, 0, 0, 0, 0, 0, 429441205967280, 0, 0, 0, 0, 0, 0, 0, 0], [−442071829672200, 0, 0, 0, 0, 0, 0, 707314927475520, 0, 0, 0, 0, 0, 0, 0], [−71573534327880, 0, 0, 0, 0, 0, 0, 0, 286294137311520, 0, 0, 0, 0, 0, 0], [−84440686566600, 0, 0, 0, 0, 0, 0, 0, 0, 135105098506560, 0, 0, 0, 0, 0], [−76425977333160, 0, 0, 0, 0, 0, 0, 0, 0, 0, 101901303110880, 0, 0, 0, 0], [−63508910741640, 0, 0, 0, 0, 0, 0, 0, 0, 0, 0, 169357095311040, 0, 0, 0], [−11301084367560, 0, 0, 0, 0, 0, 0, 0, 0, 0, 0, 0, 90408674940480, 0, 0], [−18255597824520, 0, 0, 0, 0, 0, 0, 0, 0, 0, 0, 0, 0, 24340797099360, 0], [−11834993865240, 0, 0, 0, 0, 0, 0, 0, 0, 0, 0, 0, 0, 0, 23669987730480]]  **(4.3.5)**

## ▼ 1. Given a fan, say if it is associated to P(Q) for some Q

Let us now introduce the procedures *recognizing the fan* of a wps following Section 2 in [RT].

# 1. Recognizing an admissible fan

*IsWPSFan* determines if a given list of $n+1$ vectors in $\mathbf{Z}^n$ generates a fan of a wps, by checking condition (3) in [RT, Thm.2.1].

```
> IsWPSFan := proc(V)
    local n, MM, SS, TT, co, i :
    n := nops(V[1]) :
    MM := Transpose(Matrix(V)) :
    SS := [ seq(Determinant(Matri(MM, i)), i = 1..n + 1) ] :
    TT := [seq(abs(SS[i]), i = 1..n + 1)] :
    co := 0 : for i from 1 to n + 1 do if SS[i] = 0 then co := 1 end if end do:
    if co = 0 and igcd(op(TT)) = 1 and add(TT[i]·V[i], i = 1..n + 1) = [seq(0, i = 1
        ..n)] then true else false end if:
    end proc:
```

**Examples:**

> $V := [[2, 4, 1, 2], [2, 1, 2, 2], [4, 5, 6, 1], [7, 7, 1, 1], [0, 0, 0, 1]]$

$$V := [[2, 4, 1, 2], [2, 1, 2, 2], [4, 5, 6, 1], [7, 7, 1, 1], [0, 0, 0, 1]] \qquad (5.1.1)$$

> *IsWPSFan*(*V*)

$$false \qquad (5.1.2)$$

> $V1 := [[-2,-4,-1,-2], [2, 1, 2, 2], [4, 5, 6, 1], [7, 7, 1, 1], [0, 0, 0, 1]]$

$$V1 := [[-2, -4, -1, -2], [2, 1, 2, 2], [4, 5, 6, 1], [7, 7, 1, 1], [0, 0, 0, 1]] \qquad (5.1.3)$$

> *IsWPSFan*(*V1*)

$$false \qquad (5.1.4)$$

> $V0 := [[-2,-4,-1,-2], [1, 0, 0, 0], [0, 1, 0, 0], [0, 0, 1, 0], [0, 0, 0, 1]]$

$$V0 := [[-2, -4, -1, -2], [1, 0, 0, 0], [0, 1, 0, 0], [0, 0, 1, 0], [0, 0, 0, 1]] \qquad (5.1.5)$$

> *IsWPSFan*(*V0*)

$$true \qquad (5.1.6)$$

> $V := Fan([5, 3, 2, 2, 4, 16])$
$$V := [[-3, -1, -1, -2, -5], [5, 1, 1, 2, 3], [0, 1, 0, 0, 0], [0, 0, 1, 0, 0], [0, 0, 0, 1, 0], \qquad (5.1.7)$$
$$[0, 0, 0, 0, 1]]$$

> *IsWPSFan*(*V*)

$$true \qquad (5.1.8)$$

# 2. Determining the weights vector Q

By using *IsWPSFan* the following procedure *FWeights* recognizes the fan of a wps $\mathbf{P}(Q)$ and compute the weigths vector *Q*.

```
> FWeights := proc(V)
    local n, MM, SS, TT :
    if IsWPSFan(V) = false then error "It is not a fan of a WPS"
    else
    n := nops(V[1]) :  MM := Transpose(Matrix(V)) :
    SS := [seq(Determinant(Matri(MM, i)), i = 1..n + 1)] :
    TT := [seq(abs(SS[i]), i = 1..n + 1)] :
    TT
    end if
```

```
        end proc:
```
**Examples:**
```
> V1
```
$$[[-2, -4, -1, -2], [2, 1, 2, 2], [4, 5, 6, 1], [7, 7, 1, 1], [0, 0, 0, 1]] \qquad (5.2.1)$$
```
> FWeights(V1)
Error, (in FWeights) It is not a fan of a WPS
> V0
```
$$[[-2, -4, -1, -2], [1, 0, 0, 0], [0, 1, 0, 0], [0, 0, 1, 0], [0, 0, 0, 1]] \qquad (5.2.2)$$
```
> FWeights(V0)
```
$$[1, 2, 4, 1, 2] \qquad (5.2.3)$$
```
> Q := [3, 2, 7, 181, 22]:
> V := Fan(Q);
```
$$V := [[-2, -3, -61, -8], [3, 1, 1, 1], [0, 1, 0, 0], [0, 0, 1, 0], [0, 0, 0, 1]] \qquad (5.2.4)$$
```
> FWeights(V)
```
$$[3, 2, 7, 181, 22] \qquad (5.2.5)$$

## 3. Given a weights vector Q and a fan, say if the latter is a fan of P(Q)

The following procedure *IsFanQ* answers the question if a given fan *V* defines the wps of given weights *Q*.

```
> IsFanQ := proc(V, Q)
   local Q1 :
   Q1 := FWeights(V) :
   if {op(Q)} = {op(Q1)} then true else false end if
   end proc:
```
**Examples:**
```
> Q := [3, 2, 7, 181, 22]:
> V := Fan(Q);
```
$$V := [[-2, -3, -61, -8], [3, 1, 1, 1], [0, 1, 0, 0], [0, 0, 1, 0], [0, 0, 0, 1]] \qquad (5.3.1)$$
```
> IsFanQ(V, Q)
```
$$true \qquad (5.3.2)$$
```
> Q1 := [181, 22, 2, 7, 3]
```
$$Q1 := [181, 22, 2, 7, 3] \qquad (5.3.3)$$
```
> IsFanQ(V, Q1)
```
$$true \qquad (5.3.4)$$
```
> Q2 := [3, 2, 4, 8, 1]
```
$$Q2 := [3, 2, 4, 8, 1] \qquad (5.3.5)$$
```
> IsFanQ(V, Q2)
```
$$false \qquad (5.3.6)$$

## 2. Given a polytope, say if it is associated to P(Q) for some Q

Let us then pass to the *recognition process* for polytopes defining a polarized wps $(P(Q), \mathcal{O}(m))$.

### 1. Technicalities

Given a square matrix *W* with coprime entries, the following procedure *WPP* computes a list containing
1) the list of associated pseudo-weights, which are effective weights when the columns of *W* and the origin are the vertices of a polytope defining a polarized wps,
2) the *what* of *W* as defined in [RT, Def. 3.10],
3) the associated pseudo-fan, which for general *W* may be gnerated by *rational* vectors in $Q^n$,
4) the associated polarization.

```
> WPP := proc(W1)
    local W, m, M, n, M1, i, KK, M2, IQ, qzero, deQ, M12, vzero, fan3 :
    m := igcd(op(convert(W1, list))) :
    W := (1/m) · W1 :
    n := op(W)[1] :
    M := convert(sign(Determinant(W)) · Adjoint(W), listlist) :
    M1 := [ ] : for i from 1 to n do KK := igcd(op(M[i])) :
       M1 := [op(M1), (1/KK) · M[i]] :
    end do: M12 := Matrix(M1) :
    M2 := Transpose(M12) :
    IQ := convert(M12.W, listlist) :
    qzero := abs(Determinant(M2)) :
    deQ := ilcm(qzero, seq(IQ[k, k], k = 1..n)) :
    vzero := - (1/qzero) · add((deQ/IQ[r, r]) · M1[r], r = 1..n) :
    fan3 := [vzero, op(M1)] :
    [[qzero, seq(deQ/IQ[h, h], h = 1..n)], M2, fan3, m]
  end proc:
```

The different outputs of *WPP* will be individually recalled in the following by using the following procedures:

```
> WPseudoWeights := proc(W) WPP(W)[1] end proc:
> What := proc(W) WPP(W)[2] end proc:
> PseudoFan := proc(W) WPP(W)[3] end proc:
> Polarization := proc(W) WPP(W)[4] end proc:
```

**Examples**

```
> WW := RandomMatrix(3, 3)
```

$$WW := \begin{bmatrix} 27 & 99 & 92 \\ 8 & 29 & -31 \\ 69 & 44 & 67 \end{bmatrix} \quad (6.1.1)$$

```
> WPP(WW)
```

$$\left[ [107088635536, 327244, 327244, 327244], \begin{bmatrix} -3307 & 2675 & 1649 \\ 2585 & 4539 & -5643 \\ 5737 & -1573 & 9 \end{bmatrix}, \left[ \left[ -\frac{1017}{327244}, \right. \right. \right. \quad (6.1.2)$$

$$\left. \left. -\frac{1481}{327244}, -\frac{4173}{327244} \right], [-3307, 2585, 5737], [2675, 4539, -1573], [1649, -5643,$$

```
> WW2 := QpolMat([2, 2, 7, 4, 3, 14])
```

$$WW2 := \begin{bmatrix} 42 & -6 & 0 & -14 & 0 \\ -42 & 12 & 0 & 0 & 0 \\ 0 & 0 & 21 & 0 & 0 \\ 0 & 0 & 0 & 28 & 0 \\ 0 & 0 & 0 & 0 & 6 \end{bmatrix}$$ (6.1.3)

```
> WPP(WW2)
```

$$\left[ [2, 2, 7, 4, 3, 14], \begin{bmatrix} 2 & 2 & 0 & 0 & 0 \\ 1 & 2 & 0 & 0 & 0 \\ 0 & 0 & 1 & 0 & 0 \\ 1 & 1 & 0 & 1 & 0 \\ 0 & 0 & 0 & 0 & 1 \end{bmatrix}, [[-9, -8, -2, -6, -7], [2, 1, 0, 1, 0], [2, 2, 0, 1, 0], [0, 0, 1, 0, 0], [0, 0, 0, 1, 0], [0, 0, 0, 0, 1]], 1 \right]$$ (6.1.4)

## 2. Admissible matrices

The following procedure says if a given matrix is *P*-admissible [RT, Def. 3.8] i.e. if it is the weighted transverse matrix of a wps fan. This is done by checking that the obtained pseudo-fan is actually a fan generated by integer vectors.

```
> IsAdmissible := proc(W)
    local V, r, i :
    V := PseudoFan(W)[1] :
    r := 0 : for i from 1 to nops(V) do if floor(V[i]) ≠ V[i] then r := 1 end if:
    end do:
    if r = 1 then false else true end if
  end proc:
```

### Examples

```
> WW := RandomMatrix(3, 3)
```

$$WW := \begin{bmatrix} 57 & -76 & -32 \\ 27 & -72 & -74 \\ -93 & -2 & -4 \end{bmatrix}$$ (6.2.1)

```
> IsAdmissible(WW)
```

$$false$$ (6.2.2)

```
> WW2 := QpolMat([2, 2, 7, 4, 3, 19])
```

$$WW2 := \begin{bmatrix} 798 & -114 & 0 & -266 & -42 \\ -798 & 228 & 0 & 0 & 0 \\ 0 & 0 & 399 & 0 & 0 \\ 0 & 0 & 0 & 532 & 0 \\ 0 & 0 & 0 & 0 & 84 \end{bmatrix}$$ (6.2.3)

> $IsAdmissible(WW2)$

$$true \qquad (6.2.4)$$

## 3. Recognizing an admissible polytope

The following procedure recognizes a polytope defining a polarized wps by running the previous procedure on the associated polytope matrix.

> $IsPAdmissible := \mathbf{proc}(P)$
> $\quad \mathbf{local}\ W, m:$
> $\quad W := PolytopeMatrix(P):\ \mathrm{m} := igcd(op(convert(W, list)))\ :\ IsAdmissible\left(\frac{1}{m} \cdot W\right)$
> $\quad \mathbf{end\ proc}:$

**Examples**

> $P := [[2, 4, 5], [6, 6, 1], [3, 1, 2], [5, 1, 1]]$
$$P := [[2, 4, 5], [6, 6, 1], [3, 1, 2], [5, 1, 1]] \qquad (6.3.1)$$

> $IsPAdmissible(P)$
$$false \qquad (6.3.2)$$

> $Q := [2, 2, 7, 4, 3, 19]:$
> $PP2 := Qpolytope(Q)$
$PP2 := [[0, 0, 0, 0, 0], [798, -798, 0, 0, 0], [-114, 228, 0, 0, 0], [0, 0, 399, 0, 0], [\ \ (6.3.3)$
$\quad -266, 0, 0, 532, 0], [-42, 0, 0, 0, 84]]$

> $IsPAdmissible(PP2)$
$$true \qquad (6.3.4)$$

> $mPP2 := [seq(3 \cdot PP2[i], i = 1 .. nops(PP2))]$
$mPP2 := [[0, 0, 0, 0, 0], [2394, -2394, 0, 0, 0], [-342, 684, 0, 0, 0], [0, 0, 1197, 0, 0], \ \ (6.3.5)$
$\quad [-798, 0, 0, 1596, 0], [-126, 0, 0, 0, 252]]$

> $IsPAdmissible(mPP2)$
$$true \qquad (6.3.6)$$

> $PP20 := [seq(PP2[i] + [1, 1, 1, 1, 1], i = 1 .. nops(PP2))]$
$PP20 := [[1, 1, 1, 1, 1], [799, -797, 1, 1, 1], [-113, 229, 1, 1, 1], [1, 1, 400, 1, 1], [\ \ (6.3.7)$
$\quad -265, 1, 1, 533, 1], [-41, 1, 1, 1, 85]]$

> $IsPAdmissible(PP20)$
$$true \qquad (6.3.8)$$

## 4. Determining the associated weights vector, fan and polarization

The following procedures *PPolarization, PWeights* and *PFan* recognize the polytope of a polarized wps, by using *IsAdmissible*, and in the affermative case return the polarization *m*, the weights vector *Q* and a fan of $P(Q)$, respectively. The last two outputs are obtained by running *WPP*.

```
> PPolarization := proc(P) local W, W1, m :
    W1 := PolytopeMatrix(P) : m := igcd(op(convert(W1, list))) : W := (1/m)·W1 :
    if IsAdmissible(W) = false then error "It is not an admissible polytope"
      else m end if
      end proc:
> PWeights := proc(P)
    local W, W1, m :
    W1 := PolytopeMatrix(P) : m := igcd(op(convert(W1, list))) : W := (1/m)·W1 :
    if IsAdmissible(W) = false then error "It is not an admissible polytope"
    else WPseudoWeights(W)
    end if
    end proc:

> PFan := proc(P)
    local W, W1, m :
    W1 := PolytopeMatrix(P) : m := igcd(op(convert(W1, list))) : W := (1/m)·W1 :
    if IsAdmissible(W) = false then error "It is not an admissible polytope"
    else PseudoFan(W)
    end if
    end proc:
```

**Examples**

```
> P := [[3, 5, 1], [8, 9, 2], [12, 6, 4], [4, 1, 5]]
```
$$P := [[3, 5, 1], [8, 9, 2], [12, 6, 4], [4, 1, 5]] \tag{6.4.1}$$

```
> PWeights(P)
```
Error, (in PWeights) It is not an admissible polytope

```
> Q
```
$$[2, 2, 7, 4, 3, 19] \tag{6.4.2}$$

```
> PP2 := Qpolytope(Q)
```
$PP2 := [[0, 0, 0, 0, 0], [798, -798, 0, 0, 0], [-114, 228, 0, 0, 0], [0, 0, 399, 0, 0], [-266, 0, 0, 532, 0], [-42, 0, 0, 0, 84]]$ (6.4.3)

```
> PPolarization(PP2)
```
$$1 \tag{6.4.4}$$

```
> PWeights(PP2)
```
$$[2, 2, 7, 4, 3, 19] \tag{6.4.5}$$

```
> PF := PFan(PP2)
```
$PF := [[-9, -8, -2, -6, -14], [2, 1, 0, 1, 1], [2, 2, 0, 1, 1], [0, 0, 1, 0, 0], [0, 0, 0, 1, 0], [0, 0, 0, 0, 1]]$ (6.4.6)

```
> PP20 := [seq(PP2[i] + [1, 1, 1, 1, 1], i = 1..nops(PP2))]
```
$PP20 := [[1, 1, 1, 1, 1], [799, -797, 1, 1, 1], [-113, 229, 1, 1, 1], [1, 1, 400, 1, 1], [-265, 1, 1, 533, 1], [-41, 1, 1, 1, 85]]$ (6.4.7)

```
> PWeights(PP20)
```
$$[2, 2, 7, 4, 3, 19] \tag{6.4.8}$$

```
> PPolarization(PP20)
```
$$1 \tag{6.4.9}$$

```
> PFan(PP20)
    [[-9, -8, -2, -6, -14], [2, 1, 0, 1, 1], [2, 2, 0, 1, 1], [0, 0, 1, 0, 0], [0, 0, 0, 1, 0], [0,    (6.4.10)
     0, 0, 0, 1]]
> PWeights(mPP2)
                               [2, 2, 7, 4, 3, 19]                                                    (6.4.11)
> PFan(mPP2)
    [[-9, -8, -2, -6, -14], [2, 1, 0, 1, 1], [2, 2, 0, 1, 1], [0, 0, 1, 0, 0], [0, 0, 0, 1, 0], [0,    (6.4.12)
     0, 0, 0, 1]]
> PPolarization(mPP2)
                                         3                                                            (6.4.13)
```

# 3. Equivalences

In the following two fans are said equivalent if they defines isomorphic wps's. Analogously two polytopes are said equivalent if they defines isomorphic polarized wps's, meaning that two polytopes giving different polarizations on the same wps are not considered equivalent. We will give several procedures which are able to determine the equivalence of fans and polytopes and to give concrete transformations connecting equivalent toric data.

## 1. Equivalent fans

The following procedure establishes if two fans are associated to the same reduced weights vector $Q$. First it checks if the input fans are admissible, by using *IsWPSFan*, then it compares the associated reduced weights vectors.

```
> AreEquivalentFans := proc( fan1, fan2 )
    local W1, W2 :
    if IsWPSFan( fan1 ) = false then error "The first one is not a fan of a WPS"
    else if IsWPSFan( fan2 ) = false then error "The second one is not a fan of a WPS"
    else W1 := sort(ReducedWeights(FWeights( fan1 ))) : W2
       := sort(ReducedWeights(FWeights( fan2 ))) :
    if W1 = W2  then true else false end if
    end if end if
    end proc:
```

**Examples**

```
> fan1 := [[3, 5, 1], [8, 9, 2], [12, 6, 4], [4, 1, 5]]
                   fan1 := [[3, 5, 1], [8, 9, 2], [12, 6, 4], [4, 1, 5]]                              (7.1.1)
> IsWPSFan( fan1 )
                                         false                                                        (7.1.2)
> fan2 := [[5, 2, 1], [6, 7, 5], [4, 9, 1], [8, 0, 0]]
                   fan2 := [[5, 2, 1], [6, 7, 5], [4, 9, 1], [8, 0, 0]]                               (7.1.3)
> IsWPSFan( fan2 )
                                         false                                                        (7.1.4)
> fan3 := Fan([3, 1, 5, 4])
                   fan3 := [[-1, -2, -2], [3, 1, 2], [0, 1, 0], [0, 0, 1]]                            (7.1.5)
> fan4 := Fan([7, 2, 1, 6])
                   fan4 := [[-2, -1, -2], [7, 3, 4], [0, 1, 0], [0, 0, 1]]                            (7.1.6)
> fan5 := Fan([5, 9, 3, 12])
```

$$fan5 := [[-9, -6, -6], [5, 3, 2], [0, 1, 0], [0, 0, 1]] \tag{7.1.7}$$

> $ReducedWeights(FWeights(fan5))$

$$[5, 3, 1, 4] \tag{7.1.8}$$

> $AreEquivalentFans(fan1, fan2)$

```
Error, (in AreEquivalentFans) The first one is not a fan of
a WPS
```

> $AreEquivalentFans(fan3, fan1)$

```
Error, (in AreEquivalentFans) The second one is not a fan of
a WPS
```

> $AreEquivalentFans(fan3, fan4)$

$$false \tag{7.1.9}$$

> $AreEquivalentFans(fan3, fan5)$

$$true \tag{7.1.10}$$

## Compute switching matrices connecting equivalent fans

Let us first of all reduce the fan generators to the elementary generators of the 1-skeleton. The so obtained equivalent fan is called the *reduced fan*.

> $ReducedFan := \mathbf{proc}(fan)$
> $\left[seq\left(\frac{1}{igcd(op(fan[i]))} \cdot fan[i], i = 1..nops(fan)\right)\right]$
> **end proc**:

### Examples

> $fan1$

$$[[3, 5, 1], [8, 9, 2], [12, 6, 4], [4, 1, 5]] \tag{7.1.1.1}$$

> $ReducedFan(fan1)$

$$[[3, 5, 1], [8, 9, 2], [6, 3, 2], [4, 1, 5]] \tag{7.1.1.2}$$

> $fan5$

$$[[-9, -6, -6], [5, 3, 2], [0, 1, 0], [0, 0, 1]] \tag{7.1.1.3}$$

> $fan6 := ReducedFan(fan5)$

$$fan6 := [[-3, -2, -2], [5, 3, 2], [0, 1, 0], [0, 0, 1]] \tag{7.1.1.4}$$

> $fan3$

$$[[-1, -2, -2], [3, 1, 2], [0, 1, 0], [0, 0, 1]] \tag{7.1.1.5}$$

> $AreEquivalentFans(fan5, fan6)$

$$true \tag{7.1.1.6}$$

> $AreEquivalentFans(fan3, fan6)$

$$true \tag{7.1.1.7}$$

If $S1$, $S2$ are lists of lenght k with the same underlying multiset *perm(S1,S2,\*)* is the permutation $\sigma$ of $\{1,...,k\}$ such that $S1 = [seq(S2[\sigma(i)], i = 1..n)]$

> $perm := \mathbf{proc}(S1, S2, i)$
> **local** $j, t, k$ :
> **for** $j$ **from** 1 **to** $nops(S1)$ **do if** $S2[j] = S1[i]$ **then** $t := 0$ :
> **for** $k$ **from** 1 **to** $i - 1$ **do if** $perm(S1, S2, k) = j$ **then** $t := 1$ : **end if end do**:
> **if** $t = 0$ **then** $j$ :**break**: **end if** : **end if end do**: **end proc**:

### Example

> $L1 := [1, 1, 5, 3, 3, 7, 3]$ :

```
> L2 := [3, 5, 1, 1, 3, 3, 7]:
> for i from 1 to 7 do print([i, perm(L1, L2, i)]) end do:
```
$$[1, 3]$$
$$[2, 4]$$
$$[3, 2]$$
$$[4, 1]$$
$$[5, 5]$$
$$[6, 7]$$
$$[7, 6]$$ (7.1.1.8)

```
> L3 := [seq(L2[perm(L1, L2, i)], i = 1 .. nops(L1))]
```
$$L3 := [1, 1, 5, 3, 3, 7, 3]$$ (7.1.1.9)

*PermCol(M,σ)* permutes the columns of the matrix *M* according to the permutation σ; *PermMat(n,σ)* is the *n*×*n* permutation matrix corresponding to σ.

```
> PermCol := proc(M, σ) local M1 : M1 := convert(Transpose(M), listlist) :
    Transpose(Matrix([seq(M1[σ(i)], i = 1 .. nops(M1))])) end proc:
> PermMat := proc(n, σ) PermCol(Matrix(n, shape = identity), σ) end proc:
```

**Examples**

```
> σ := proc(i) if i = 1 then 3 else if i = 2 then 4 else if i = 3 then 2 else if i = 4 then 4 fi
    fi fi fi end proc:
> M := RandomMatrix(2, 4)
```
$$M := \begin{bmatrix} 45 & -38 & 87 & -98 \\ -81 & -18 & 33 & -77 \end{bmatrix}$$ (7.1.1.10)

```
> PermCol(M, σ)
```
$$\begin{bmatrix} 87 & -98 & -38 & -98 \\ 33 & -77 & -18 & -77 \end{bmatrix}$$ (7.1.1.11)

```
> PermMat(4, σ)
```
$$\begin{bmatrix} 0 & 0 & 0 & 0 \\ 0 & 0 & 1 & 0 \\ 1 & 0 & 0 & 0 \\ 0 & 1 & 0 & 1 \end{bmatrix}$$ (7.1.1.12)

Given two equivalent fans *fan1* and *fan2*, let *RedFan1* and *RedFan2* their reductions. The procedure *FSwitchMatrices* computes matrices α, β, γ, δ such that
- α and β are diagonal matrices with positive integer entries such that
  *FanMatrix(fan1) = α · FanMatrix(RedFan1)*
  and *FanMatrix(fan2) = β · FanMatrix(RedFan2)*.
- γ is an invertible integer matrix and δ is a permutation matrix such that
  *FanMatrix(RedFan2) = γ · FanMatrix(RedFan1) · δ*, according with [RT, Prop.2.7].

```
> FSwitchMatrices := proc(fan1, fan2)
    local K1, K2, α, β, γ, δ, W1, W2, RW1, RW2, F, RedFan1, RedFan2, V1, V2, Vzero1,
      Vzero2:
```

```
        if AreEquivalentFans(fan1, fan2) = false then
            error "The two fans are not equivalent"
        else K1 := [seq(igcd(op(fan1[i])), i = 1..nops(fan1))] : K2
                := [seq(igcd(op(fan2[i])), i = 1..nops(fan2))] :
        α := Matrix(Vector(K1), shape = diagonal) : β := Matrix(Vector(K2), shape
               = diagonal) :
        W1 := FWeights(fan1) : W2 := FWeights(fan2) :
        RW1 := ReducedWeights(W1); RW2 := ReducedWeights(W2) :
        F := proc(i) perm(RW1, RW2, i) end proc:
        δ := PermCol(Matrix(nops(fan1), shape = identity), F) :
        RedFan1 := ReducedFan(fan1) : RedFan2 := ReducedFan(fan2) :
        V1 := FanMatrix(RedFan1) : V2 := FanMatrix(RedFan2) :
        Vzero1 := Vzero(V1) : Vzero2 := Vzero(V2.δ) :
        γ := Vzero1.MatrixInverse(Vzero2) :
        [α, β, γ, δ]
        end if:
        end proc:
```

# Examples

> $B := RandomMatrix(7)$

$$B := \begin{bmatrix} 70 & 42 & -67 & -95 & 65 & 31 & -2 \\ -32 & 18 & 22 & -20 & 86 & -50 & 50 \\ -1 & -59 & 14 & -25 & 20 & -80 & 10 \\ 52 & 12 & 16 & 51 & -61 & 43 & -16 \\ -13 & -62 & 9 & 76 & -48 & 25 & -9 \\ 82 & -33 & 99 & -44 & 77 & 94 & -50 \\ 72 & -68 & 60 & 24 & 9 & 12 & -22 \end{bmatrix}$$ (7.1.1.13)

> $C := HermiteForm(B, output = 'U')$;

$C := [\,[59523579227, 97498373629, 210296252447, 234258218556,$ (7.1.1.14)
    $151884931336, 97390657747, -264296715696],$
    $[1871447117, 3065391104, 6611805279, 7365179874, 4775328039,$
    $3062004470, -8309603237],$
    $[48093966353, 78776907600, 169915536342, 189276367911,$
    $122720254260, 78689875132, -213546925722],$
    $[49139979620, 80490255380, 173611091497, 193393008874,$
    $125389341961, 80401330011, -218191435925],$
    $[56152476216, 91976577642, 198386176826, 220991062982,$
    $143282966266, 91874962222, -249328337355],$
    $[2246013516, 3678923005, 7935144887, 8839305901, 5731100399,$
    $3674858543, -9972753712],$
    $[67456847740, 110492900968, 238324416405, 265480019619,$
    $172128067905, 110370828775, -299522030436]\,]$

> $Determinant(C)$

$$-1 \tag{7.1.1.15}$$

```
> Q := [6, 8, 6, 7, 7, 9, 25, 2]
```
$$Q := [6, 8, 6, 7, 7, 9, 25, 2] \tag{7.1.1.16}$$

```
> Q1 := [9, 25, 8, 2, 7, 6, 7, 6]
```
$$Q1 := [9, 25, 8, 2, 7, 6, 7, 6] \tag{7.1.1.17}$$

```
> fan7 := Fan(Q)
```
$fan7 := [[-29, -1, -6, -6, -16, -9, -10], [6, 0, 1, 1, 3, 1, 2], [0, 1, 0, 0, 0, 0, 0],$ (7.1.1.18)
$[18, 0, 4, 3, 9, 3, 6], [0, 0, 0, 1, 0, 0, 0], [0, 0, 0, 0, 1, 0, 0], [0, 0, 0, 0, 0, 1, 0], [0, 0, 0, 0, 0, 0, 1]]$

```
> fan8 := Fan(Q1)
```
$fan8 := [[-25, -12, -3, -23, -9, -23, -9], [9, 4, 1, 8, 3, 8, 3], [0, 1, 0, 0, 0, 0, 0],$ (7.1.1.19)
$[0, 0, 1, 0, 0, 0, 0], [0, 0, 0, 1, 0, 0, 0], [0, 0, 0, 0, 1, 0, 0], [0, 0, 0, 0, 0, 1, 0], [0, 0, 0, 0, 0, 0, 1]]$

```
> F8 := MatrixFan(C.FanMatrix(fan8))
```
$F8 := [[-9905176817293, -311423050140, -8003202205882,$ (7.1.1.20)
$-8177266778230, -9344199607473, -373753751016, -11225333104614],$
$[3451957617350, 108531042906, 2789118793877, 2849780358785,$
$3256456861703, 130253314164, 3912032549496], [97498373629,$
$3065391104, 78776907600, 80490255380, 91976577642, 3678923005,$
$110492900968], [210296252447, 6611805279, 169915536342,$
$173611091497, 198386176826, 7935144887, 238324416405],$
$[234258218556, 7365179874, 189276367911, 193393008874,$
$220991062982, 8839305901, 265480019619], [151884931336, 4775328039,$
$122720254260, 125389341961, 143282966266, 5731100399, 172128067905]$
$, [97390657747, 3062004470, 78689875132, 80401330011, 91874962222,$
$3674858543, 110370828775], [-264296715696, -8309603237,$
$-213546925722, -218191435925, -249328337355, -9972753712,$
$-299522030436]]$

```
> FSwitchMatrices(fan7, F8)
```
$$\left[\left[\begin{array}{cccccccc} 1 & 0 & 0 & 0 & 0 & 0 & 0 & 0 \\ 0 & 1 & 0 & 0 & 0 & 0 & 0 & 0 \\ 0 & 0 & 1 & 0 & 0 & 0 & 0 & 0 \\ 0 & 0 & 0 & 1 & 0 & 0 & 0 & 0 \\ 0 & 0 & 0 & 0 & 1 & 0 & 0 & 0 \\ 0 & 0 & 0 & 0 & 0 & 1 & 0 & 0 \\ 0 & 0 & 0 & 0 & 0 & 0 & 1 & 0 \\ 0 & 0 & 0 & 0 & 0 & 0 & 0 & 1 \end{array}\right], \left[\begin{array}{cccccccc} 1 & 0 & 0 & 0 & 0 & 0 & 0 & 0 \\ 0 & 1 & 0 & 0 & 0 & 0 & 0 & 0 \\ 0 & 0 & 1 & 0 & 0 & 0 & 0 & 0 \\ 0 & 0 & 0 & 1 & 0 & 0 & 0 & 0 \\ 0 & 0 & 0 & 0 & 1 & 0 & 0 & 0 \\ 0 & 0 & 0 & 0 & 0 & 1 & 0 & 0 \\ 0 & 0 & 0 & 0 & 0 & 0 & 1 & 0 \\ 0 & 0 & 0 & 0 & 0 & 0 & 0 & 1 \end{array}\right], \tag{7.1.1.21}\right.$$

$$\left[\begin{matrix} 491 & 1744 & 762 & -551 & 225 & -235 & -803 \\ 85 & -6 & 51 & -52 & 57 & -23 & -120 \\ 114 & 354 & 166 & -82 & 0 & -61 & -167 \\ 144 & 309 & 249 & -177 & 138 & -3 & -299 \\ 300 & 986 & 334 & -437 & 217 & -75 & -390 \\ 241 & 612 & -11 & -551 & 335 & 82 & -102 \\ 167 & 543 & 265 & -234 & 111 & -99 & -270 \end{matrix}\right], \left[\begin{matrix} 0 & 0 & 0 & 0 & 0 & 1 & 0 & 0 \\ 0 & 0 & 0 & 0 & 0 & 0 & 1 & 0 \\ 0 & 1 & 0 & 0 & 0 & 0 & 0 & 0 \\ 0 & 0 & 0 & 0 & 0 & 0 & 0 & 1 \\ 0 & 0 & 0 & 1 & 0 & 0 & 0 & 0 \\ 1 & 0 & 0 & 0 & 0 & 0 & 0 & 0 \\ 0 & 0 & 0 & 0 & 1 & 0 & 0 & 0 \\ 0 & 0 & 1 & 0 & 0 & 0 & 0 & 0 \end{matrix}\right]$$

> *fan77* := [[-29, -1, -6, -6, -16, -9, -10], [12, 0, 2, 2, 6, 2, 4], [0, 1, 0, 0, 0, 0, 0], [18, 0, 4, 3, 9, 3, 6], [0, 0, 0, 1, 0, 0, 0], [0, 0, 0, 0, 1, 0, 0], [0, 0, 0, 0, 0, 1, 0], [0, 0, 0, 0, 0, 0, 1]]

*fan77* := [[-29, -1, -6, -6, -16, -9, -10], [12, 0, 2, 2, 6, 2, 4], [0, 1, 0, 0, 0, 0, 0], [18, 0, 4, 3, 9, 3, 6], [0, 0, 0, 1, 0, 0, 0], [0, 0, 0, 0, 1, 0, 0], [0, 0, 0, 0, 0, 1, 0], [0, 0, 0, 0, 0, 0, 1]]     **(7.1.1.22)**

> *fan88* := [seq(i·F8[i], i = 1..8)]

*fan88* := [[-9905176817293, -311423050140, -8003202205882,     **(7.1.1.23)**
-8177266778230, -9344199607473, -373753751016, -11225333104614],
[6903915234700, 217062085812, 5578237587754, 5699560717570,
6512913723406, 260506628328, 7824065098992], [292495120887,
9196173312, 236330722800, 241470766140, 275929732926, 11036769015,
331478702904], [841185009788, 26447221116, 679662145368,
694444365988, 793544707304, 31740579548, 953297665620],
[1171291092780, 36825899370, 946381839555, 966965044370,
1104955314910, 44196529505, 1327400098095], [911309588016,
28651968234, 736321525560, 752336051766, 859697797596, 34386602394,
1032768407430], [681734604229, 21434031290, 550829125924,
562809310077, 643124735554, 25724009801, 772595801425], [
-2114373725568, -66476825896, -1708375405776, -1745531487400,
-1994626698840, -79782029696, -2396176243488]]

> *FSwitchMatrices*(*fan7*, *fan8*)

$$\left[\left[\begin{matrix} 1 & 0 & 0 & 0 & 0 & 0 & 0 & 0 \\ 0 & 1 & 0 & 0 & 0 & 0 & 0 & 0 \\ 0 & 0 & 1 & 0 & 0 & 0 & 0 & 0 \\ 0 & 0 & 0 & 1 & 0 & 0 & 0 & 0 \\ 0 & 0 & 0 & 0 & 1 & 0 & 0 & 0 \\ 0 & 0 & 0 & 0 & 0 & 1 & 0 & 0 \\ 0 & 0 & 0 & 0 & 0 & 0 & 1 & 0 \\ 0 & 0 & 0 & 0 & 0 & 0 & 0 & 1 \end{matrix}\right], \left[\begin{matrix} 1 & 0 & 0 & 0 & 0 & 0 & 0 & 0 \\ 0 & 1 & 0 & 0 & 0 & 0 & 0 & 0 \\ 0 & 0 & 1 & 0 & 0 & 0 & 0 & 0 \\ 0 & 0 & 0 & 1 & 0 & 0 & 0 & 0 \\ 0 & 0 & 0 & 0 & 1 & 0 & 0 & 0 \\ 0 & 0 & 0 & 0 & 0 & 1 & 0 & 0 \\ 0 & 0 & 0 & 0 & 0 & 0 & 1 & 0 \\ 0 & 0 & 0 & 0 & 0 & 0 & 0 & 1 \end{matrix}\right],\right.$$

**(7.1.1.24)**

$$\begin{bmatrix} -9 & 6 & 0 & 18 & -29 & 0 & 0 \\ 0 & 0 & 0 & 0 & -1 & 0 & 1 \\ -2 & 1 & 0 & 4 & -6 & 0 & 0 \\ -2 & 1 & 0 & 3 & -6 & 1 & 0 \\ -4 & 3 & 0 & 9 & -16 & 0 & 0 \\ 0 & 1 & 0 & 3 & -9 & 0 & 0 \\ -3 & 2 & 1 & 6 & -10 & 0 & 0 \end{bmatrix}, \begin{bmatrix} 0 & 0 & 0 & 0 & 0 & 1 & 0 & 0 \\ 0 & 0 & 0 & 0 & 0 & 0 & 1 & 0 \\ 0 & 1 & 0 & 0 & 0 & 0 & 0 & 0 \\ 0 & 0 & 0 & 0 & 0 & 0 & 0 & 1 \\ 0 & 0 & 0 & 1 & 0 & 0 & 0 & 0 \\ 1 & 0 & 0 & 0 & 0 & 0 & 0 & 0 \\ 0 & 0 & 0 & 0 & 1 & 0 & 0 & 0 \\ 0 & 0 & 1 & 0 & 0 & 0 & 0 & 0 \end{bmatrix}$$

> *fan7*

$[[-29, -1, -6, -6, -16, -9, -10], [6, 0, 1, 1, 3, 1, 2], [0, 1, 0, 0, 0, 0, 0], [18, 0,$ **(7.1.1.25)**
$4, 3, 9, 3, 6], [0, 0, 0, 1, 0, 0, 0], [0, 0, 0, 0, 1, 0, 0], [0, 0, 0, 0, 0, 1, 0], [0, 0,$
$0, 0, 0, 0, 1]]$

> *IsWPSFan(fan7)*

$$true \qquad (7.1.1.26)$$

> *IsWPSFan(fan77)*

$$false \qquad (7.1.1.27)$$

## 2. Equivalent polytopes

The following procedure establishes if two polytopes are associated with isomorphic polarized wps's. First it checks if the input polytopes are admissible, by using *IsPadmissible*, then it compares the associated weights vectors and the associated polarizations.

> *AreEquivalentPols* := **proc**(*pol1, pol2*)
> local *W1, W2* :
> **if** *IsPAdmissible(pol1)* = *false* **then error** "The first one is not admissible" :
> **else if** *IsPadmissible(pol2)* = *false* **then error** "The second one is not admissible" :
> **else** *W1* := *sort(PWeights(pol1))* : *W2* := *sort(PWeights(pol2))* :
> **if** *W1* = *W2* **and** *PPolarization(pol1)* = *PPolarization(pol2)* **then** *true* **else** *false* **end if**
> **end if end if end proc*:

**Examples**

> *Q*

$$[6, 8, 6, 7, 7, 9, 25, 2] \qquad (7.2.1)$$

> *P1* := *Qpolytope(Q)*;

$P1 := [[0, 0, 0, 0, 0, 0, 0], [1050, 0, -4725, 0, 0, 0, 0], [0, 2100, 0, 0, 0, 0, 0], [-300, 0,$ **(7.2.2)**
$1800, 0, 0, 0, 0], [-300, 0, 0, 1800, 0, 0, 0], [-700, 0, 0, 0, 1400, 0, 0], [-84, 0, 0,$
$0, 0, 504, 0], [-2100, 0, 0, 0, 0, 0, 6300]]$

> *Q1*

$$[9, 25, 8, 2, 7, 6, 7, 6] \qquad (7.2.3)$$

> *P2* := *Qpolytope(Q1)*;

$P2 := [[0, 0, 0, 0, 0, 0, 0], [56, 0, 0, 0, 0, 0, 0], [-700, 1575, 0, 0, 0, 0, 0], [-700, 0,$ **(7.2.4)**
$6300, 0, 0, 0, 0], [-1600, 0, 0, 1800, 0, 0, 0], [-700, 0, 0, 0, 2100, 0, 0], [-1600, 0,$
$0, 0, 0, 1800, 0], [-700, 0, 0, 0, 0, 0, 2100]]$

> $AreEquivalentPols(P1, P2)$

$$true \quad (7.2.5)$$

> $fan3$

$$[[-1, -2, -2], [3, 1, 2], [0, 1, 0], [0, 0, 1]] \quad (7.2.6)$$

> $fan5$

$$[[-9, -6, -6], [5, 3, 2], [0, 1, 0], [0, 0, 1]] \quad (7.2.7)$$

> $P3 := Fpolytope(fan3)$

$$P3 := [[0, 0, 0], [20, 0, 0], [-4, 12, 0], [-10, 0, 15]] \quad (7.2.8)$$

> $P5 := Fpolytope(fan5)$

$$P5 := [[0, 0, 0], [4, 0, 0], [-36, 60, 0], [-6, 0, 15]] \quad (7.2.9)$$

> $AreEquivalentPols(P3, P5)$

$$true \quad (7.2.10)$$

> $P30 := [seq(P3[i] + [1, 1, 2], i = 1..nops(P3))]$

$$P30 := [[1, 1, 2], [21, 1, 2], [-3, 13, 2], [-9, 1, 17]] \quad (7.2.11)$$

> $AreEquivalentPols(P5, P30)$

$$true \quad (7.2.12)$$

> $mP3 := [seq(3 \cdot P3[i], i = 1..nops(P3))]$

$$mP3 := [[0, 0, 0], [60, 0, 0], [-12, 36, 0], [-30, 0, 45]] \quad (7.2.13)$$

> $AreEquivalentPols(P5, mP3)$

$$false \quad (7.2.14)$$

## Compute switching data connecting equivalent polytopes

Given two equivalent polytopes *pol1* and *pol2*, let *CP1* and *CP2* be the $n \times n+1$ matrices whose columns are given by vertices of *pol1* and *pol2*, respectively. The procedure P*SwitchMatrices* computes matrices θ, δ and a vector *v* such that

- θ is an invertible $n \times n$ matrix (it is the transverse of the matrix α in the procedure *FSwitchMatrices*);
- δ is an $(n + 1) \times (n + 1)$ permutation matrix (it is the same matrix δ given by the procedure F*SwitchMatrices*);
- *v* is a vector in $\mathbf{Z}^n$,
- and finally $\theta \cdot CP2 \cdot \delta + v = CP1$, according with [RT, Prop.3.21].

> $PSwitchMatrices := \mathbf{proc}(pol1, pol2)$
  $\mathbf{local}$ F1, F2, δ, CP1, CP2, CCP1, pol3, CP3, θ, v :
  $\mathbf{if}$ $AreEquivalentPols(pol1, pol2) = false$ $\mathbf{then}$
      $\mathbf{error}$ "The two polytopes are not equivalent"
  $\mathbf{else}$
  $F1 := PFan(pol1) : F2 := PFan(pol2) :$
  $\theta := MatrixInverse(Transpose(FSwitchMatrices(F1, F2)[3])) :$
  $\delta := FSwitchMatrices(F1, F2)[4] :$
  $CP1 := Transpose(Matrix(pol1)) :$
  $CP2 := Transpose(Matrix(pol2)) :$
  $CCP1 := \theta.CP2.\delta :$
  $v := pol1[1] - convert(Column(CCP1, 1), list) :$
  $pol3 := [seq(convert(Transpose(CCP1), listlist)[i] + v, i = 1..nops(F1))] :$
  $CP3 := Transpose(Matrix(pol3)) :$
  $[\theta, \delta, v]$

```
        end if:
    end proc:
Example
```

> $PP3 := Fpolytope(fan3)$
$$PP3 := [[0, 0, 0], [20, 0, 0], [-4, 12, 0], [-10, 0, 15]] \qquad (7.2.1.1)$$

> $PP5 := Fpolytope(fan5)$
$$PP5 := [[0, 0, 0], [4, 0, 0], [-36, 60, 0], [-6, 0, 15]] \qquad (7.2.1.2)$$

> $P2$
$$[[0, 0, 0, 0, 0, 0, 0], [56, 0, 0, 0, 0, 0, 0], [-700, 1575, 0, 0, 0, 0, 0], [-700, 0, \qquad (7.2.1.3)$$
$$6300, 0, 0, 0, 0], [-1600, 0, 0, 1800, 0, 0, 0], [-700, 0, 0, 0, 2100, 0, 0], [$$
$$-1600, 0, 0, 0, 0, 1800, 0], [-700, 0, 0, 0, 0, 0, 2100]]$$

>
> $PSM := PSwitchMatrices(PP3, PP5)$

$$PSM := \left[ \begin{bmatrix} 1 & 1 & 0 \\ -3 & -2 & -2 \\ 0 & 0 & 1 \end{bmatrix}, \begin{bmatrix} 0 & 0 & 1 & 0 \\ 1 & 0 & 0 & 0 \\ 0 & 1 & 0 & 0 \\ 0 & 0 & 0 & 1 \end{bmatrix}, [-4, 12, 0] \right] \qquad (7.2.1.4)$$

> $CP1 := Transpose(Matrix(PP3)); CP2 := Transpose(Matrix(PP5));$

$$CP1 := \begin{bmatrix} 0 & 20 & -4 & -10 \\ 0 & 0 & 12 & 0 \\ 0 & 0 & 0 & 15 \end{bmatrix}$$

$$CP2 := \begin{bmatrix} 0 & 4 & -36 & -6 \\ 0 & 0 & 60 & 0 \\ 0 & 0 & 0 & 15 \end{bmatrix} \qquad (7.2.1.5)$$

> $P3 := [seq(convert(Transpose(PSM[1].CP2.PSM[2]), listlist)[i] + PSM[3], i = 1$
$\quad ...nops(PP3))]:$
> $evalb(P3 = PP3)$
$$true \qquad (7.2.1.6)$$

# ▼ 4. Cohomology

Finally we want end up this collection of procedures on wps, by computing the cohomology of line bundles and sheaves of holomorphic forms. It is a well known fact that such a computation for toric varieties reduces to a count of lattice points in a polytope.

## ▼ 1. Lattice points in a polytope

Given a polytope *pol*, the following procedure *PPP* computes a list of two vectors: *Lm* is the vector of the smallest coordinates of the points in *pol*, *LM* is the vector of the greatest coordinates in *pol*.

> $PPP := \mathbf{proc}(pol)$
$\quad \mathbf{local}\ n, m, M, Lm, LM, i:$
$\quad n := nops(pol[1]):$

```
    m := proc(i) min(seq(pol[j][i], j = 1 .. n + 1)) end proc:
    M := proc(i) max(seq(pol[j][i], j = 1 .. n + 1)) end proc:
    Lm := [ ]: for i from 1 to n do Lm := [op(Lm), m(i)] : end do:
    LM := [ ]: for i from 1 to n do LM := [op(LM), M(i)] : end do:
    [Lm, LM]
    end proc:
```

### Example

```
> Pol := [[2, 3, 5], [-1, 4, 6], [3, 2, 5], [4, 2, 8]] :
> PPP(Pol)
```
$$[[-1, 2, 5], [4, 4, 8]] \qquad (8.1.1)$$

The following procedure determines the polytope, called the *reduced polytope*, whose vertices are obtained by reducing the vertices of a given polytope.

```
> RedPol := proc(pol)
    local rr : rr := igcd(seq(igcd(op(pol[i])), i = 1 .. nops(pol))) : [seq( (1/rr) · pol[i], i = 1 .. nops(pol) )] : end proc:

> RedPol([[4, 6, 10], [-4, 8, 12]])
```
$$[[2, 3, 5], [-2, 4, 6]] \qquad (8.1.2)$$

The following procedures determines points which either belongs or are internal to a a given polytope. Namely *IsInPolytope(point,pol)* is *true* if the integer *point* belongs to the polytope *pol*, *false* otherwise, while *IsInternal(point,pol)* verifies if the *point* is internal to the polytope *pol*:

```
> IsInPolytope := proc(point, pol)
    local P1, F, X, L, i, co, A, AA, M, n, newpoint, newpol :
    n := nops(pol[1]) :
    M := PolytopeMatrix(pol) :
    newpol := [[seq(0, t = 1 .. n)], seq(pol[k] − pol[1], k = 2 .. n + 1)] :
    newpoint := point − pol[1] :
    A := convert(sign(Determinant(M)).Adjoint(M), listlist) :
    AA := [−add(A[r], r = 1 .. nops(A)), op(A)] :
    X := Vector(newpoint);
    L := [DotProduct(X − Vector(newpol[2]), Vector(AA[1]))] :
    for i from 2 to nops(AA) do L := [op(L), DotProduct(X, Vector(AA[i]))] : end do:
    co := 0 : for i from 1 to nops(L) do if L[i] < 0 then co := 1 : end if end do:
    if co = 0 then true else false
    end if
    end proc:

> IsInternal := proc(point, pol)
    local P1, F, X, L, i, co, A, AA, M, n, newpoint, newpol :
    n := nops(pol[1]) :
    M := PolytopeMatrix(pol) :
    newpol := [[seq(0, t = 1 .. n)], seq(pol[k] − pol[1], k = 2 .. n + 1)] :
    newpoint := point − pol[1] :
    A := convert(sign(Determinant(M)).Adjoint(M), listlist) :
    AA := [−add(A[r], r = 1 .. nops(A)), op(A)] :
    X := Vector(newpoint);
    L := [DotProduct(X − Vector(newpol[2]), Vector(AA[1]))] :
    for i from 2 to nops(AA) do L := [op(L), DotProduct(X, Vector(AA[i]))] : end do:
    co := 0 : for i from 1 to nops(L) do if L[i] ≤ 0 then co := 1 : end if end do:
```

```
            if co = 0 then true else false
            end if
          end proc:
```

**Examples**

```
> Q := [2, 5, 7, 5, 2]:
> Pol := Qpolytope(Q);
     Pol := [[0, 0, 0, 0], [7, 0, 0, 0], [-5, 10, 0, 0], [-7, 0, 14, 0], [0, 0, 0, 35]]    (8.1.3)
> IsInPolytope([0, 0, 0, 0], Pol)
                                  true                                                      (8.1.4)
> IsInternal([0, 0, 0, 0], Pol)
                                  false                                                     (8.1.5)
> IsInPolytope([0, 1, 0, 0], Pol)
                                  true                                                      (8.1.6)
> IsInternal([0, 1, 0, 0], Pol)
                                  false                                                     (8.1.7)
> IsInPolytope([0, 1, 1, 1], Pol)
                                  true                                                      (8.1.8)
> IsInternal([0, 1, 1, 1], Pol)
                                  true                                                      (8.1.9)
> IsInPolytope([0, -1, 0, 10], Pol)
                                  false                                                     (8.1.10)
> Pol := [[0, 0, 0], [2, 0, 0], [0, 5, 0], [0, 0, 12]]
                Pol := [[0, 0, 0], [2, 0, 0], [0, 5, 0], [0, 0, 12]]                        (8.1.11)
> IsInPolytope([1, 1, 1], Pol)
                                  true                                                      (8.1.12)
> IsInternal([0, 1, 0], Pol)
                                  false                                                     (8.1.13)
```

The procedure *PointsInPol* computes the number of points belonging to the polytope and a list of such points, by using procedures *PPP* and *IsInPolytope*; the procedure *InternalPoints* does the same for internal points by finally applying *IsInternal*.

```
> PointsInPol := proc(pol)
    local n, Lm, LM, Lm1, LM1, Count, P, T, Q, u, L, uu, uuu :
    L := [ ]:
    n := nops(pol[1]) :
    Lm := PPP(pol)[1] :
    LM := PPP(pol)[2] :
    Lm1 := [seq(ceil(Lm[i]), i = 1..n)] :
    LM1 := [seq(ceil(LM[i]), i = 1..n)] :
    Count := 0 :
    P := Lm1 :
    T := proc(x) evalb(x > 0) end proc:
    while P ≠ LM1
      do
    if IsInPolytope(P, pol) then Count := Count + 1 : L := [op(L), P] : end if:
    Q := LM1 − P :
    if Q ≠ [seq(0, i = 1..n)] then
    uu := ListTools[SelectLast](T, Q) :
```

```
            uuu := [ListTools[SearchAll](uu, Q)] :
            u := uuu[nops(uuu)] : end if:
            P := [seq(P[i], i = 1 .. u − 1), P[u] + 1, seq(Lm1[i], i = u + 1 .. n)] :
            end do:
            if IsInPolytope(LM1, pol) then Count := Count + 1 : L := [op(L), LM1] : end if:
            [Count, L]
            end proc:

>   PointsInPol := proc(pol)
        local n, Lm, LM, Lm1, LM1, Count, P, T, Q, u, L, uu, uuu, L1, co, point, M, newpol, A,
            AA, newpoint, i, X :
        n := nops(pol[1]) :
        M := PolytopeMatrix(pol) :
        newpol := [[seq(0, t = 1 .. n)], seq(pol[k] − pol[1], k = 2 .. n + 1)] :
        A := convert(sign(Determinant(M)).Adjoint(M), listlist) :
        AA := [−add(A[r], r = 1 .. nops(A)), op(A)] :
        Lm := PPP(pol)[1] :
        LM := PPP(pol)[2] :
        Lm1 := [seq(ceil(Lm[i]), i = 1 .. n)] :
        LM1 := [seq(ceil(LM[i]), i = 1 .. n)] :
        L := [ ] :
        Count := 0 :
        T := proc(x) evalb(x > 0) end proc:
        P := Lm1 :
        while P ≠ LM1 do
        newpoint := P − pol[1] :
        X := Vector(newpoint);
        L1 := [DotProduct(X − Vector(newpol[2]), Vector(AA[1]))] :
        for i from 2 to nops(AA) do L1 := [op(L1), DotProduct(X, Vector(AA[i]))] : end do:
        co := 0 : for i from 1 to nops(L1) do if L1[i] < 0 then co := 1 : end if end do:
        if co = 0 then Count := Count + 1 : L := [op(L), P] : end if:
        Q := LM1 − P :
        if Q ≠ [seq(0, i = 1 .. n)] then
        uu := ListTools[SelectLast](T, Q) :
        uuu := [ListTools[SearchAll](uu, Q)] :
        u := uuu[nops(uuu)] : end if:
        P := [seq(P[i], i = 1 .. u − 1), P[u] + 1, seq(Lm1[i], i = u + 1 .. n)] :
        end do:
        P := LM1 :
        newpoint := P − pol[1] :
        X := Vector(newpoint);
        L1 := [DotProduct(X − Vector(newpol[2]), Vector(AA[1]))] :
        for i from 2 to nops(AA) do L1 := [op(L1), DotProduct(X, Vector(AA[i]))] : end do:
        co := 0 : for i from 1 to nops(L1) do if L1[i] < 0 then co := 1 : end if end do:
        if co = 0  then Count := Count + 1 : L := [op(L), LM1] : end if:
        [Count, L]
        end proc:

>   InternalPoints := proc(pol)
        local n, Lm, LM, Lm1, LM1, Count, P, T, Q, u, L, uu, uuu, L1, co, point, M, newpol, A,
            AA, newpoint, i, X :
```

```
n := nops(pol[1]):
M := PolytopeMatrix(pol):
newpol := [[seq(0, t = 1..n)], seq(pol[k] − pol[1], k = 2..n + 1)]:
A := convert(sign(Determinant(M)).Adjoint(M), listlist):
AA := [−add(A[r], r = 1..nops(A)), op(A)]:
Lm := PPP(pol)[1]:
LM := PPP(pol)[2]:
Lm1 := [seq(ceil(Lm[i]), i = 1..n)]:
LM1 := [seq(ceil(LM[i]), i = 1..n)]:
L := []:
Count := 0:
T := proc(x) evalb(x > 0) end proc:
P := Lm1:
while P ≠ LM1 do
newpoint := P − pol[1]:
X := Vector(newpoint);
L1 := [DotProduct(X − Vector(newpol[2]), Vector(AA[1]))]:
for i from 2 to nops(AA) do L1 := [op(L1), DotProduct(X, Vector(AA[i]))]: end do:
co := 0: for i from 1 to nops(L1) do if L1[i] ≤ 0 then co := 1: end if end do:
if co = 0 then Count := Count + 1: L := [op(L), P]: end if:
Q := LM1 − P:
if Q ≠ [seq(0, i = 1..n)] then
uu := ListTools[SelectLast](T, Q):
uuu := [ListTools[SearchAll](uu, Q)]:
u := uuu[nops(uuu)]: end if:
P := [seq(P[i], i = 1..u − 1), P[u] + 1, seq(Lm1[i], i = u + 1..n)]:
end do:
P := LM1:
newpoint := P − pol[1]:
X := Vector(newpoint);
L1 := [DotProduct(X − Vector(newpol[2]), Vector(AA[1]))]:
for i from 2 to nops(AA) do L1 := [op(L1), DotProduct(X, Vector(AA[i]))]: end do:
co := 0: for i from 1 to nops(L1) do if L1[i] ≤ 0 then co := 1: end if end do:
if co = 0 then Count := Count + 1: L := [op(L), LM1]: end if:
[Count, L]
end proc:
```

## Examples

> $Q := [2, 3, 5]$

$$Q := [2, 3, 5] \tag{8.1.14}$$

> $Pol := Qpolytope(Q)$

$$Pol := [[0, 0], [5, 0], [-3, 6]] \tag{8.1.15}$$

> $Pol := [[0, 0], [0, 3], [5, 0]]$

$$Pol := [[0, 0], [0, 3], [5, 0]] \tag{8.1.16}$$

> $PointsInPol(Pol)$

$$[13, [[0, 0], [0, 1], [0, 2], [0, 3], [1, 0], [1, 1], [1, 2], [2, 0], [2, 1], [3, 0], [3, 1], [4, 0], [5, 0]]] \tag{8.1.17}$$

> $InternalPoints(Pol)$

$$[4, [[1, 1], [1, 2], [2, 1], [3, 1]]] \tag{8.1.18}$$

The following procedure $H(point, pol)$ returns the dimension of the smallest face of *pol* containing the given *point*. This the key point to provide a procedure computing Bott-Tu formulas. This corresponds precisely to the function $s(u, D)$ in [RT, Theorem 4.8].

```
> H := proc(point, pol)
    local n, P1, F, X, L, i, co, M, newpol, newpoint, A, AA:
    n := nops(pol[1]):
    if not IsInPolytope(point, pol) then 0:
    else if IsInternal(point, pol) then n:
    else
    newpol := [[seq(0, t = 1..n)], seq(pol[k] - pol[1], k = 2..n + 1)]:
    newpoint := point - pol[1]:
    M := Matri(Transpose(Matrix(newpol)), 1):
    A := convert(sign(Determinant(M)).Adjoint(M), listlist):
    AA := [-add(A[r], r = 1..nops(A)), op(A)]:
    P1 := Vector(pol[1]):
    F := PFan(pol):
    X := Vector(point);
    L := [DotProduct(X - Vector(pol[2]), Vector(AA[1]))]:
    for i from 2 to nops(F) do L := [op(L), DotProduct(X, Vector(AA[i]))]: end do:
    co := 0: for i from 1 to nops(L) do if L[i] = 0 then co := co + 1: end if end do:
    n - co
    end if end if:
    end proc:

> H([1, 12/5], Pol)
```
$$1 \tag{8.1.19}$$

```
> H([1, 1], Pol)
```
$$2 \tag{8.1.20}$$

```
> H([5, 0], Pol)
```
$$0 \tag{8.1.21}$$

```
> H([12, 0], Pol)
```
$$0 \tag{8.1.22}$$

## 2. Cohomology of line bundles

The following procedure $hO(q,m,Q)$ computes the complex dimension of complex $q$--th cohomogy vector space $H^q(P(Q), \mathcal{O}(m))$ where $\mathcal{O}(1)$ is the generator of Pic$(P(Q))$, and $m$ is an integer. This procedure produces a result for any $q \geq 0$ and for any $m \geq 0$. Moreover if $P(Q)$ is Gorenstein then $hO(q, m, Q)$ produces a result also for $m \leq -\frac{|Q|}{\delta}$, where $|Q| := \sum_{i=0}^{n} q_i$ and $\delta := lcm(Q)$ (computed by procedures *Qsum* and *Qdelta* in the Preliminaries). It is obtained by Proposition 4.6 in [RT]. Although it admits large ranges of uncomputability on the polarization $m$, it is actually faster than the following powerful procedure $hOmega(q,p,m,Q)$.

```
> hO := proc(q, m, Q) local h; if q < 0 then error "Cohomology of negative level" end if:
    if m ≥ 0 then if q = 0 then h
        := PointsInPol(MatrixPolytope(HermiteForm(PolytopeMatrix(Qpolytope(Q,
        m)))))[1] else h := 0 end if: elif IsGorenstein(Q) = false
        then print("Unable to compute it: no Gorenstein case, please use hOmega") elif m ≤
```

$$-\frac{Qsum(Q)}{Qdelta(Q)} \text{ then if } q = nops(Q) - 1 \text{ then } h := hO\left(0, -m - \frac{Qsum(Q)}{Qdelta(Q)}, Q\right)$$

**else** $h := 0$ **end if:** **else** $print\left(-\frac{Qsum(Q)}{Qdelta(Q)} < m < 0\right)$ :

$print($"Unable to compute it in this range. Please use hOmega" $)$ **end if: end proc**:

**Examples**

> $Q$

$$[2, 3, 5] \tag{8.2.1}$$

> $hO(-1, 1, [1, 1, 1])$
Error, (in hO) Cohomology of negative level

> $hO(0, 1, Q)$

$$21 \tag{8.2.2}$$

> $hO(1, 1, Q)$

$$0 \tag{8.2.3}$$

> $hO(0, -1, Q)$

"Unable to compute it: no Gorenstein case, please use hOmega" (8.2.4)

> $hO(0, -1, [1, 1, 1])$

*true*

"Unable to compute it in this range. Please use hOmega" (8.2.5)

> $hO(2, -1, [1, 1, 1])$

*true*

"Unable to compute it in this range. Please use hOmega" (8.2.6)

> $hO(0, 0, [1, 1, 1])$

$$1 \tag{8.2.7}$$

> $hO(2, -3, [1, 1, 1])$

$$1 \tag{8.2.8}$$

> $hO(1, -1, [1, 1, 1])$

*true*

"Unable to compute it in this range. Please use hOmega" (8.2.9)

# 3. Bott-Tu formulas

The following procedure *hOmega(q,p,m,Q)* computes the complex dimension of complex *q*--th cohomogy vector space $H^q(P(Q), \Omega^p(m))$ where *m* is a twisting integer. This procedure produces a result for any $q \geq 0, p \geq 0$ and for any *m*. It is obtained by Theorem 4.9 in [RT]. Although it is a more powerful procedure than *hO(q,m,Q)* it is unfortunately slower.

> $hOmega := \text{proc}(q, p, m, Q)$ **local** $h, PP, i$; **if** $q < 0$ **then**
  **error** "Cohomology of negative level"**end if:** **if** $p < 0$ **then**
  **error** "Cohomology of negative level differential forms"**end if:** **if** $m = 0$ **then if** $q = p$
  **then** $h := 1$ **else** $h := 0$ **end if:** **elif** $m > 0$ **then if** $q = 0$ **then** $h := 0 : PP$
  $:= PointsInPol(MatrixPolytope(HermiteForm(PolytopeMatrix(Qpolytope(Q, m)))))$ : **for** $i$ **to** $PP[1]$ **do** $h := h + \text{binomial}(H(PP[2][i], Qpolytope(Q, m)), p)$
  **end do:else** $h := 0$ **end if:** **elif** $m < 0$ **then if** $q = nops(Q) - 1$ **then** $h := 0 : PP$
  $:= PointsInPol(Qpolytope(Q, -m))$ : **for** $i$ **to** $PP[1]$ **do** $h := h$
  $+ \text{binomial}(H(PP[2][i], Qpolytope(Q, -m)), nops(Q) - 1 - p)$ **end do:else** $h$
  $:= 0$ **end if: end if:** $h$; **end proc:**

**Examples**

```
> hOmega(-1, 0, 1, [1, 1, 1])
Error, (in hOmega) Cohomology of negative level
> hOmega(0, 0, 2, Q)
```
$$71 \tag{8.3.1}$$

```
> hO(0, 2, Q)
```
$$71 \tag{8.3.2}$$

```
> hO(0, 1, Q)
```
$$21 \tag{8.3.3}$$

```
> hOmega(0, 0, 1, Q)
```
$$21 \tag{8.3.4}$$

```
> hO(0, 3, Q)
```
$$151 \tag{8.3.5}$$

```
> hOmega(0, 0, 3, Q)
```
$$151 \tag{8.3.6}$$

```
> hOmega(1, 0, 1, Q)
```
$$0 \tag{8.3.7}$$

```
> hOmega(0, 0, -1, Q)
```
$$0 \tag{8.3.8}$$

```
> hOmega(0, 0, -1, [1, 1, 1])
```
$$0 \tag{8.3.9}$$

```
> hOmega(2, 0, -1, [1, 1, 1])
```
$$0 \tag{8.3.10}$$

```
> hOmega(0, 0, 0, [1, 1, 1])
```
$$1 \tag{8.3.11}$$

```
> hOmega(2, 0, -3, [1, 1, 1])
```
$$1 \tag{8.3.12}$$

```
> hOmega(1, 0, -1, [1, 1, 1])
```
$$0 \tag{8.3.13}$$

```
> hOmega(2, 0, -1, Q)
```
$$11 \tag{8.3.14}$$

```
> hO(2, -1, Q)
```
$$\text{"Unable to compute it: no Gorenstein case, please use hOmega"} \tag{8.3.15}$$

```
> hOmega(0, 1, 1, [1, 2, 2, 3])
```
$$13 \tag{8.3.16}$$

```
> hOmega(0, 2, 1, [1, 2, 2, 3])
```
$$3 \tag{8.3.17}$$

```
> hO(0, 1, [1, 1, 2])
```
$$4 \tag{8.3.18}$$

```
> hO(0, 3, [1, 1, 2])
```
$$16 \tag{8.3.19}$$

```
> hOmega(0, 2, 3, [1, 1, 2])
```
$$4 \tag{8.3.20}$$

```
> hO(0, 1, [1, 3, 3, 2])
```
$$11 \tag{8.3.21}$$

> $hO(0, 2, [1, 3, 3, 2])$

42     (8.3.22)

> $hOmega(0, -1, 1, [1, 1, 1, 1])$
Error, (in hOmega) Cohomology of negative level differential forms

# References

[RT] M. Rossi, L. Terracini *Weighted projective spaces from the toric point of view with computational applications,*
http://www2.dm.unito.it/paginepersonali/rossi/wps.pdf
[Cox] http://www.cs.amherst.edu/~dac/toric.html